\documentclass{amsart}
\usepackage{amsmath,amssymb,latexsym}

\def\acl{\mathrm{acl}}

\def\M{\mathfrak M}
\def\N{\mathfrak N}

\def\ale{\lesssim}

\def\len{\unlhd}

\def\C{\tilde C}
\def\Z{\tilde Z}
\def\N{\tilde N}
\def\H{\mathcal H}

\def\ZZ{\mathbb Z}
\def\SU{\mathrm{SU}}

\def\Ind#1#2{#1\setbox0=\hbox{$#1x$}\kern\wd0\hbox to 0pt{\hss$#1\mid$\hss}
\lower.9\ht0\hbox to 0pt{\hss$#1\smile$\hss}\kern\wd0}

\def\Notind#1#2{#1\setbox0=\hbox{$#1x$}\kern\wd0\hbox to 0pt{\mathchardef
\nn="3236\hss$#1\nn$\kern1.4\wd0\hss}\hbox to 0pt{\hss$#1\mid$\hss}\lower.9\ht0
\hbox to 0pt{\hss$#1\smile$\hss}\kern\wd0}

\theoremstyle{plain}
\newtheorem{theorem}{Theorem}[section]
\newtheorem{proposition}[theorem]{Proposition}
\newtheorem{fact}[theorem]{Fact}
\newtheorem{lemma}[theorem]{Lemma}
\newtheorem{corollary}[theorem]{Corollary}

\theoremstyle{definition}
\newtheorem{definition}[theorem]{Definition}
\newtheorem{remark}[theorem]{Remark}
\newtheorem{expl}[theorem]{Example}

\def\bsp{\begin{expl}}
\def\ebsp{\end{expl}}
\def\beh{\begin{claim}}
\def\ebeh{\end{claim}}
\def\defn{\begin{definition}}
\def\edefn{\end{definition}}
\def\satz{\begin{theorem}}
\def\esatz{\end{theorem}}
\def\tats{\begin{fact}}
\def\etats{\end{fact}}
\def\kor{\begin{corollary}}
\def\ekor{\end{corollary}}
\def\lmm{\begin{lemma}}
\def\elmm{\end{lemma}}
\def\bem{\begin{remark}}
\def\ebem{\end{remark}}
\def\bew{\par\noindent{\em Proof: }}

\def\satzli{\begin{proposition}}
\def\esatzli{\end{proposition}}

\title{A Fitting Theorem for Simple Theories}
\author{Daniel Palac\'in \and Frank O. Wagner}
\address{Universit\"at M\"unster; Institut f\"ur Mathematische Logik und
Grundlagenforschung,
Einsteinstrasse 62, 48149 M\"unster, Germany}
\email{daniel.palacin@uni-muenster.de}
\address{Universit\'e de Lyon; CNRS; Universit\'e Lyon 1; Institut Camille
Jordan UMR5208, 43 bd du 11 novembre 1918, 69622 Villeurbanne Cedex, France}
\email{wagner@math.univ-lyon1.fr}
\thanks{The first author was supported by the project SFB 878 and the
project MTM 2011-26840. The second author was partially supported by ValCoMo
(ANR-13-BS01-0006)}

\begin{document}
\begin{abstract} The Fitting subgroup of a type-definable group in
a simple theory is relatively definable and nilpotent. Moreover, the Fitting
subgroup of a supersimple hyperdefinable group has a normal hyperdefinable
nilpotent subgroup of bounded index, and is itself of bounded index in a
hyperdefinable subgroup.\end{abstract}
\maketitle

\section{Introduction}
The {\em Fitting subgroup} $F(G)$ of a group $G$ is the group generated by all
normal nilpotent subgroups. Since the product of two normal nilpotent subgroups
of class $c$ and $c'$ respectively is again a normal nilpotent subgroup of
class $c+c'$, it is clear that the Fitting subgroup of a finite group is
nilpotent. In general, this need not be the case, and some additional
finiteness conditions are needed. For groups with the chain condition on
centralisers ($\M_c$), nilpotency of the Fitting subgroup was shown by
Bryant \cite{br79} for periodic groups, by Poizat and Wagner
\cite{PW93,wa95} in the stable case (stability being a model-theoretic
condition implying, in particular, that definable groups are $\M_c$),
and generally by Derakhshan and Wagner \cite{DW97}.

In this paper, we shall prove nilpotency of the Fitting subgroup in the
model-theoretic context of groups type-definable in simple structures (Theorem
\ref{fitting}), and virtual nilpotency of the Fitting subgroup for groups
hyperdefinable in supersimple structures (Theorem \ref{hyperfitting}).
Simplicity is a fundamental broadening of stability mentioned above; algebraic
examples
of simple structures include $\omega$-free bounded pseudo-algebraically closed
fields, generic (differential) difference fields, or vector spaces over a finite
field with a non-degenerate bilinear form. 

In this context, a weaker
chain condition on centralisers one might call (uniform) $\widetilde{\M}_c$
arises naturally \cite[Theorem 4.2.12]{wa00}: The (uniform) chain condition on
centralisers up to finite index. More precisely, we shall assume that there are
natural numbers $n,d<\omega$ such that any chain $(C_G(a_j:j<i):i<k)$ of
centralizers, each of index at least $d$ in its predecessor, has length $k<n$.

Similarly to the approach in \cite{PW93,wa95} we shall also need this chain
condition on quotients by relatively definable subgroups, which again follows
from simplicity, as the quotients are again type-definable.
In this context, in contrast to $\M_c$-groups, the centralizer of an infinite
set need not be definable, and the centralizer of a relatively definable
set (or even subgroup) need not even be type-definable. Nevertheless, the
general theory of simplicity (Fact \ref{RelDef}) yields the (relative)
definability of a weaker, still purely group theoretic notion introducted by
Haimo \cite{ha53} under the name of FC-centralizer: For $H\le G$ put
$$\C_G(H)=\{g\in G:|H:C_H(g)|<\infty\}.$$
However, whereas for subgroups $H,K\le G$ trivially $H\le
C_G(K)$ if and only if $K\le C_G(H)$, no such symmetry has to hold for the
FC-centralisers, even if one only asks for inclusion up to finite index.
Nevertheless, for type-definable groups in a simple theory, symmetry does hold
(Proposition \ref{sym}). Note that Fact \ref{RelDef} (and more generally Fact
\ref{dcc}) and symmetry are the only consequences of simplicity of the ambient
theory we will use until we study hyperdefinable groups in section four.
In fact, Hempel \cite{he} has recently shown that
symmetry of the FC-centraliser holds in general for type-definable groups;
adapting our methods she deduces nilpotency of the Fitting
subgroup of any group with the uniform $\widetilde\M_c$-condition in definable
sections.

In the course of the paper, we shall also need that soluble and nilpotent
groups are contained (up to finite index) in relatively definable soluble and
nilpotent supergroups. In the definable simple case this has been shown by
Milliet \cite{mi}. We follow his proofs (which adapt
ideas from the hyperdefinable case \cite{wa00}) in the type-definable context,
adding some details about the existence of a suitable normal relatively
definable abelian/central series.

Some words about definability: A subset of a structure is {\em definable} if it
is the set of realizations of some first-order formula, which may contain
parameters, such as
$$\begin{aligned} Z(G)&=\{x:G(x)\land\forall y\,[G(y)\to
xy=yx]\},\quad\text{or}\\
C_G(a)&=\{x:G(x)\land xa=ax\},\end{aligned}$$
where the group $G$ is given by a predicate $G(x)$, or by a formula
abbreviated as $G(x)$. A set is definable {\em relative} to a superset $X$ if
it is the intersection with $X$ of a definable set; it is {\em
type-definable} if it is given as the intersection of definable sets. Finally,
$X$ is {\em hyperdefinable} if it is the quotient of a type-definable set in
countably many variables by a type-definable equivalence relation.
In particular, a type-definable group $(G,\cdot)$ is a type-definable 
set $G$ together with a definable subset of $G\times G\times G$ which is the 
graph of the group operation. A subgroup $H$ of a group $G$ is 
relatively definable if the underlying set $H$ is definable relative to $G$.
Note that the quotient of a type-definable group by a relatively definable
normal subgroup is still type-definable\footnote{or {\em type-interpretable},
if we do not consider classes modulo $\emptyset$-definable equivalence relations
as elements}; however the quotient of a type-definable group by a type-definable
normal subgroup is only hyperdefinable, which in a way makes hyperdefinable
groups a more natural category to work with. 

Recall that if $\kappa$ is an infinite cardinal at least the size of the
language, we call a structure $\kappa$-saturated if any intersection of less
than
$\kappa$ definable sets is non-empty as soon as all finite subintersections are.
In the presence of type-definable sets, we shall always assume that the
ambient structure is more highly saturated than the size of the intersections
used in the type-definitions.

Although we make no assumption on the language used, we can in fact assume that
it is reduced to the pure group language $\{\cdot,{}^{-1},=\}$ (in addition to 
the predicate(s) (type-)defining the ambient group). All our
definable sets have an easy group-theoretic interpretation (except for Section
four where the model theory becomes more involved).

\section{Almost just definitions}\label{sec1}
\defn A subgroup $H$ of $G$ is {\em almost contained} in another subgroup $K$,
denoted by $H\ale K$, if $H\cap K$ has finite index in $H$. Then $H$ and $K$
are {\em commensurable}, denoted by $H\sim K$, if $H\ale K$ and $K\ale H$.\edefn
Observe that $\ale$ is a transitive relation among subgroups of $G$, and that
$\sim$ is an equivalence relation.
\defn Let $H$ be a subgroup of a group $G$.
The {\em almost normalizer} of $H$ is defined as $$\N_G(H)=\{g\in G: H\sim
H^g\}.$$\edefn
Note that if $H$ and $K$ are commensurable, then $\N_G(H)=\N_G(K)$. By
compactness, if $G$ is type-definable and $H$ relatively definable, then
$\N_G(H)$ is relatively definable if and only if all $\N_G(H)$-conjugates of
$H$ are uniformly commensurable with $H$; otherwise $\N_G(H)$ is not even
type-definable.

We shall usually work in a context where commensurability is uniform. Then a
theorem by Schlichting \cite{sch}, generalized by Bergmann and Lenstra
\cite{BL} (see also \cite[Theorem 4.2.4]{wa00} for definability questions),
provides an invariant object:
\tats\label{sch} Let $\H$ be
a family of uniformly commensurable subgroups of a group $G$, i.e.\ the index
$|H:H\cap H^*|$ is finite and bounded independently of $H,H^*\in\H$. Then there
is a subgroup $N$ commensurable with any $H\in\H$, which is invariant under any
automorphism of $G$ stabilizing $\H$ setwise. Moreover, $N$ is a finite
extension of a finite intersection of groups in $\H$; if the latter are
relatively definable, so is $N$.\etats

It follows that if $H$ is uniformly commensurable with all its
$\N_G(H)$-conjugates, then there is a normal subgroup $\tilde H\len\N_G(H)$
commensurable with $H$. Clearly, any two choices for $\tilde H$ will be
commensurable, and $N_G(\tilde H)=\N_G(\tilde H)=\N_G(H)$.
We call $\tilde H$ a
{\em conjugacy-connected component} of $H$. Clearly, any two
conjugacy-connected components of $H$ are commensurable, and have the same
normalizer, which is equal to their almost normalizer.

\defn\label{Adefn} Let $K$ and $H$ be subgroups of a group $G$ with $H\ale
\N_G(K)$, and suppose $K$ has a conjugacy-connected component $\tilde K$.
The {\em almost centralizer of $H$ modulo $K$} is given by
$$\C_G(H/K)=\{g\in \N_G(K):|H:C_H(g/\tilde K|\text{ is finite}\}.$$
For $n<\omega$ the {\em $n$-th iterated almost centralizer of $H$ modulo $K$} is
defined inductively by $\C_G^0(H/K)=K$,
and if $\C_G^n(H/K)$ has a
conjugacy-connected component, then
$$\C_G^{n+1}(H/K)=\bigcap_{i\le
n}\N_G(\C_G^i(H/K))\cap\C_G(H/\C_G^n(H/K)).$$
If $K=\{1\}$ it is omitted.
Moreover, we define the {\em almost centre} of $G$ as $\Z(G)=\C_G(G)$ and for
$n<\omega$ the {\em $n$-th iterated almost centre of $G$} as
$\Z_n(G)=\C_G^n(G)$. \edefn
Thus $\C_G(H/K)=\C^1_G(H/K)$ and $\Z(G)=\C_G^1(G)$.
For any subgroup $L$ we put $\C^n_L(H/K)=\C^n_G(H/K)\cap L$ and
$\N_L(K)=\N_G(K)\cap L$.

\bem One easily sees that $\C_G(H/K)$ is an $H$-invariant subgroup of
$\N_G(K)$, as are all $\C^n(H/K)$ for $n>0$. Moreover, if $H^*\sim H$ and
$K^*\sim K$, then $\C^n_G(H/K)=\C^n_G(H^*/K^*)$ for all $n>0$.\ebem

\tats\label{RelDef} If $H$ is type-definable and $K$ a relatively definable
subgroup of $H$ in a simple theory, then commensurability among
conjugates is uniform. In particular, $K$ has a conjugacy-connected
component, and both $\N_G(K)$ and $\C^n_G(H/K)$ are relatively
definable subgroups of $G$.\etats
\bew This follows immediately from \cite[Lemma 4.2.6]{wa00}.\qed

In a simple theory, the existence of generic elements yields the following
symmetry property, which plays an essential role throughout the paper.
\satzli\label{sym} Let $G$ be type-definable in a simple theory, and $H$
and $K$ be type-definable subgroups of $G$. The following are
equivalent:\begin{enumerate}
\item $H\ale\C_G(K)$.
\item There are independent generic elements $h\in H$ and $k\in K$ with
$[h,k]=1$.
\end{enumerate}
In particular, $H\ale\C_G(K)$ if and only if $K\ale\C_G(H)$.\esatzli

\bew Suppose $H\ale\C_G(K)$. So there is a generic $h\in H$ with $h\in\C_G(K)$.
Thus $C_K(h)$ has finite index in $K$, and there is generic $k\in K$ over
$h$ with $k\in C_G(h)$. Then $h$ and $k$ are independent, and $[h,k]=1$.

Conversely, suppose $h\in H$ and $k\in K$ are independent generics with
$[h,k]=1$. As $k\in C_K(h)$ and $k$ is generic over $h$, the index of $C_K(h)$
in $K$ is finite. Thus $h\in\C_G(K)$; as $h$ is generic, we get
$H\ale\C_G(K)$.\qed

Of course, if $H,K\le N_G(N)$ we also have $H\ale\C_G(K/N)$ if and only if
$K\ale\C_G(H/N)$, by working in the group $N_G(N)/N$. Thus symmetry also holds
for relative almost centralizers.

We shall finish this section by recalling two group-theoretic facts.
\tats[{\cite[Theorem 3.1]{neu}}]\label{neum} There is a finite bound
of the size of
conjugacy classes in a group $G$ if and only if the derived subgroup $G'$ is
finite.\etats
\tats[\cite{B52,ros}]\label{baer} Let $H$ and $N$ be subgroups of $G$ with $N$
normalized by $H$. If the set of commutators
$$\{[h,n]:h\in H,n\in N\}$$
is finite, then the group $[H,N]$ is finite.\etats

\section{Nilpotency in type-definable groups in a simple theory}
In this section, we shall generalize the results of Milliet \cite{mi} to
the relatively definable context. For this we need the following result.
\tats[{\cite[Theorem 4.2.12]{wa00}}]\label{dcc} Let $G$ be
type-definable in a simple theory, and $\H$ a family of uniformly relatively
definable subgroups. Then there are $n,d<\omega$ such that any
chain $(\bigcap_{j<i}H_j:i<k)$ of intersections of groups $H_i\in\H$, each of
index at least $d$ in its predecessor, has length $k<n$.\etats

\lmm\label{soluble} Let $G$ be a type-definable group in a simple theory, and
$H$ a soluble subgroup of $G$. Then there is a relatively definable soluble
subgroup $S$ containing $H$, and a series of relatively definable
$S$-invariant subgroups
$$\{1\}=S_0<S_1<\cdots<S_n=S,$$
all normalized by $N_G(H)$, such that $S_i/S_{i+1}$ is abelian for all
$i<n$. The derived length $n$ of $S$ is at most three times the derived length
of $H$. Moreover, $S_1$ and $S_n/S_{n-1}$ are finite.\elmm
\bew Suppose first that $H$ is abelian. By Fact \ref{dcc} there is a finite
tuple $\bar h\in H$ and $d<\omega$ such
that for any $h\in H$
$$|C_G(\bar h):C_G(\bar h,h)|<d.$$
Hence the $N_G(H)$-conjugates of $C_G(\bar h)$ are all
commensurable, and $N_G(H)\le \N_G(C_G(\bar h))=:N$. Now Fact \ref{RelDef}
yields that $N$ is relatively definable, and $C_G(\bar h)$ has a relatively
definable conjugacy-connected component $C$ normalized by $N$. Then
$$H=C_H(\bar h)\le C_G(\bar h)\ale C.$$
Now, as
$$C\sim C_G(\bar h)\sim C_G(\bar h,h)=C_{C_G(\bar h)}(h)\sim C_C(h)$$
for any
$h\in H$, the relatively definable subgroup $\Z(C)$ of $G$
contains $H\cap C$, whence $H\ale \Z(C)$. By compactness there is a finite
bound on the size of conjugacy classes in
$\Z(C)$, so $\Z(C)'$ is finite by Fact \ref{neum} and hence definable.
Put $S_2=C_{\Z(C)}(\Z(C)')$, a relatively definable
characteristic subgroup of
$C$, which must be normalized by $N$. Since $\Z(C)'$ is finite,
$S_2$ has finite index in $\Z(C)$, so $H\ale S_2$. If
$S_1=S_2\cap\Z(C)'$, then $S_1$ is finite, abelian
(even central in $S_2$)
and normalized by $N$, and
$S_2/S_1$ is abelian. Put $S_3=HS_2$, a finite extension of $S_2$ and thus
relatively definable. Then $S_3/S_2$ is abelian as well;
moreover, if $\bar h'$ is a system of representatives of $S_3/S_2$, then
$$N_G(H)\le N_1=\{g\in N:h^g\in S_3 \mbox{ for all $h\in\bar h'$} \}\le
N_G(S_3).$$
So $N_1$ is a relatively definable subgroup of $G$ normalizing $S_3$ and
containing $N_G(H)$. Replace $G$ by $N_1/S_3$, we finish by
induction.\qed

Note that the above proof merely uses the $\widetilde{\mathfrak M}_c$-condition
for $G$ and for certain relatively definable sections of $G$, but not symmetry
of the almost centraliser. This is different for nilpotency, where the
following lemma is used.
\lmm\label{CZG} Let $G$ be a type-definable group in a simple theory. Then there
is a characteristic relatively definable subgroup $G_0$ of finite index and a
finite characteristic subgroup $N\le Z(G_0)$ such that $\Z(G)\le
C_G(G_0/N)$.\elmm
\bew As trivially $\C_G(G)\ale\C_G(G)$, Proposition \ref{sym} yields
$$G\ale \C_G(\C_G(G))=\C_G(\Z(G)),$$ and so $\C_G(\Z(G))$ is a
characteristic relatively definable subgroup of finite index in $G$. For
independent $g\in\Z(G)$ and $h\in\C_G(\Z(G))$ we have
$$[g,h]\in\acl(g)\cap\acl(h)=\acl(\emptyset).$$
As every element in $\Z(G)$ is the product of two generic elements $g,g'$
each of which can be chosen independently of $h\in\C_G(\Z(G))$, and
$$[gg',h]=[g,h][[g,h],g'][g',h]\in\acl(\emptyset),$$
the set of commutators
$$\{[g,h]:g\in\Z(G),h\in\C_G(\Z(G))]\}$$
is bounded, whence finite by compactness. By Fact \ref{baer} the
characteristic group $Z=[\Z(G),\C_G(\Z(G))]$ is finite. We put
$G_0=\C_G(\Z(G))\cap C_G(Z)$ and $N=Z\cap G_0$.\qed

\lmm\label{nilpotent} Let $G$ be a type-definable group in a simple theory, and
$H$ a nilpotent subgroup of $G$. Then there is a relatively definable nilpotent
subgroup $N$ with $H\ale N$, and a series of relatively definable
$N_G(H)$-invariant subgroups $$\{1\}=N_0<N_1<\cdots<N_n=N$$
such that $N_{i+1}\le C_G(N/N_i)$ for all $i<n$. The nilpotency class $n$ of $N$
is at most two times the nilpotency class of $H$.\par
The same conclusion holds if $H$ is merely $FC$-nilpotent, i.e.\ $H\ale\Z_k(H)$
for some $k$.\elmm
\bew We use induction on the (FC-)nilpotency class of $H$.
Consider the family $\H=\{C_G(g):g\in \C_G(H)\}$. By Fact
\ref{dcc} there is a finite intersection $C$ of groups in $\H$ such that any
further intersection has boundedly finite index. As $\H$ is clearly
$N_G(H)$-invariant, we obtain $\N_G(C)\ge N_G(H)$. Fact \ref{RelDef} yields
that $\N_G(C)$ is relatively definable, and $C$ has a relatively definable
conjugacy-connected component $\tilde C$. So $N_G(\tilde C)=\N_G(C)\ge N_G(H)$.
As $H\ale C_G(g)$ for any $C_G(g)\in\H$, we get
$$H\ale C\sim\tilde C.$$
By construction
$$\C_G(H\cap\tilde C)=\C_G(H)\le\C_G(C)=\C_G(\tilde C),$$
whence $\Z(H\cap\tilde C)\le\Z(\tilde C)$. Now consider the characteristic
relatively definable groups $G_0$ and $N$ given by Lemma
\ref{CZG} applied to $\tilde C$. Then $H\ale G_0$. We put $N_1=N$ and
$N_2=\Z(\tilde C)\cap G_0$. Note that $N_1$, $N_2$ and $G_0$ are all relatively
definable and normalized by $N_G(\tilde C)\ge N_G(H)$, and $G_0\cap H$ is
normalized by $N_G(H)$; moreover $N_1\le C_G(G_0\cap H)$
and $N_2\le C_G(G_0\cap H/N_1)$. Now $(H\cap G_0)/N_2$ is a nilpotent subgroup
of $N_G(\tilde C)/N_2$ of smaller (FC-)nilpotency class, and we finish by
induction.\qed

\bem\label{nilpotentnormal} If in addition $H$ is
normalized by $N$, then $NH$ is nilpotent of class at most
three times the class of $H$; if $c$ is the nilpotency class of $H$ and $\bar
h\in H$ is a system of representatives of $NH/N$, then
$$\begin{aligned}\{1\}&\le C_{N_1}(\bar h)\le C_{N_1}^2(\bar h)\le\cdots\le
C_{N_1}^{c+1}(\bar h)=N_1\\
&\le C_{N_2}(\bar h/N_1)\le C_{N_2}^2(\bar h/N_1)\le\cdots\le
C_{N_2}^{c+1}(\bar h)=N_3\\
&\le\cdots\le C_{N_n}^{c+1}(\bar h/N_{n-1})=N_n=N\\
&\le NZ(H)\le NZ_2(H)\le\cdots\le NZ_c(H)=NH\end{aligned}$$
is a relatively definable central series for $NH$ normalised by $N_G(H)$.\ebem

\lmm\label{FGsoluble} Let $G$ be a
type-definable group in a simple theory and $\H$ a directed system of nilpotent
subgroups. Then $H=\bigcup\H$ is soluble.\elmm
\bew Let $(n,d)$ be the bounds given by Fact \ref{dcc} for chains
of centralizers, and $(n',d')$ the bounds for chains of centralizers modulo
$\Z(G)$. The proof is by induction on $n$.

Consider $G_0$ and $N$ as given by Lemma \ref{CZG}. As $HG_0/G_0$ is
finite and nilpotent, we may assume that $H\le G_0$. If $H\le\Z(G)$ we are done,
as $(\Z(G)\cap G_0)/N$ and $N$ are abelian.

If $H\not\le\Z(G)$, then consider some $H_0\in\H$ with $H_0\not\le\Z(G)$ and we
take $h_0\in C_{H_0}(H_0/\Z(G))\setminus\Z(G)$; note that such an  element
exists
as $H_0/\Z(G)$ is nilpotent. If $C_H(h_0/\Z(G))$ has index
greater than $d'$ in $H$, then there is some $H_1>H_0$ in $\H$ such that
$C_{H_1}(h_0/\Z(G))$ has index greater than $d'$ in $H_1$. Choose $h_1\in
C_{H_1}(H_1/\Z(G))\setminus\Z(G)$, and note that $C_G(h_1,h_0/\Z(G))$ has
index
greater than $d'$ in $C_G(h_1/\Z(G))$. If $C_H(h_1/\Z(G))$ has index more than
$d'$ in $H$, then we can iterate this process, which must stabilize after at
most $n'$ steps. It follows that there is some $h\in
H\setminus\Z(G)$ such that $C_H(h/\Z(G))$ has index at most $d'$ in $H$.

Since $C_G(h)$ has infinite index in $G$, the induction hypothesis for $n-1$
yields that $C_H(h)$ is soluble. Moreover, as $N$ is central in $G_0$ the
map from $C_H(h/N)$ to $N$ given by $x\mapsto[h,x]$ is a homomorphism with
abelian image and kernel $C_H(h)$. Thus $C_H(h/N)/C_H(h)$ is abelian. Similarly,
as $\Z(G)$ is centralised by $H$ modulo $N$, the map $x\mapsto[h,x]N$ is a
homomorphism from $C_H(h/\Z(G))$ to $\Z(G)/N$ with abelian image and kernel
$C_H(h/N)$. Therefore, $C_H(h/\Z(G))$ is soluble. Finally, as
$C_H(h/\Z(G))$ contains a normal subgroup $K$ of $H$ with $H/K$ finite,
whence nilpotent, we see that $H$ must be soluble.\qed

\kor A locally nilpotent subgroup $H$ of a type-definable group in a simple
theory is soluble.\ekor
\bew The collection of finitely generated subgroups of $H$ satisfies the
hypotheses of Lemma \ref{FGsoluble}.\qed

\lmm\label{3} Let $G$ be a type-definable group acting definably on a
type-definable
abelian group $A$, in a simple theory. Suppose that $H\le G$ is abelian,
and that there are $\bar g=(g_i:i<k)$ in $H$
and $m_i<\omega$ for $i<k$ such that $(g_i-1)^{m_i}A$ is finite for all
$i<k$, and for any $g\in H$ the index of $C_A(\bar g,g)$ in $C_A(\bar g)$ is
finite. Put $m=1+\sum_{i<k}(m_i-1)$. Then there is a relatively
definable supergroup $\bar H$ of $H$ such that $\C_A^m(\bar H)$ has finite index
in $A$.\elmm
\bew By \cite[Lemma 4.2.6]{wa00} the group
$$\bar H=\{g\in C_G(\bar g):|C_A(\bar g):C_A(\bar g,g)| \ \mbox{finite}\}$$
is relatively definable, and it clearly contains $H$. By the pigeonhole
principle, for any $m$ indices
$(i_j:j<m)\in k^m$, there must be at least one
$i<k$ such that $m_i$ of the indices are equal to
$i$. As the group ring $\ZZ(H)$ is commutative, this implies that
$$(g_{i_0}-1)(g_{i_1}-1)\cdots(g_{i_{m-1}}-1)A$$
is finite. Hence there is a subgroup $A_0$ of finite index in $A$ such
that
$$(g_{i_0}-1)(g_{i_1}-1)\cdots(g_{i_{m-1}}-1)A_0=0$$
for all
of the finitely many choices
$(i_j:j<m)\in k^m$. It follows that for all choices of
$(i_j:1\le j<m)\in k^{m-1}$ we have
$$(g_{i_1}-1)\cdots(g_{i_{m-1}}-1)A_0\le C_A(g_i:i<k).$$
As $C_A(h_0,g_i:i<k)$ has finite index in $C_A(g_i:i<k)$ for all $h_0\in
\bar H$, the group
$$(h_0-1)(g_{i_1}-1)\cdots(g_{i_{m-1}}-1)A_0$$
is finite for all choices of $(i_j:1\le j<m)\in k^{m-1}$, as is
$$(h_0-1)(g_{i_1}-1)\cdots(g_{i_{m-1}}-1)A.$$
By the same argument
(keeping $h_0$ fixed) and the fact that $\bar H\le
C_G(\bar g)$ we see that for any $h_1$ in $G$ and all choices of $(i_j:2\le
j<m)\in k^{m-2}$ the group
$$(h_1-1)(h_0-1)(g_{i_2}-1)\cdots(g_{i_{m-1}}-1)A$$
is finite, and inductively that
$$(h_{m-1}-1)\cdots(h_1-1)(h_0-1)A$$
is finite for any $(h_j:j<m)$ in $\bar H$.
It follows that
$$\bar H\le\C_{\bar H}((h_{m-2}-1)\cdots(h_0-1)A)$$
for all $(h_j:j<m-1)$ in $G$, whence by symmetry
$$(h_{m-2}-1)\cdots(h_0-1)A\ale\C_A(\bar H).$$
But $\C_A(\bar H)$ is relatively definable; we may divide out and note that
$$(h_{m-2}-1)\cdots(h_0-1)A/\C_A(\bar H)$$
is finite for all choices of $(h_j:j<m-1)$ in $\bar H$. Hence
$$\bar H\le\C_{\bar H}((h_{m-3}-1)\cdots(h_0-1)A/\C_A(\bar H))$$
and by symmetry
$$(h_{m-3}-1)\cdots(h_0-1)A\ale\C_A(\bar H/\C_A(\bar H))=\C_A^2(\bar H).$$
Inductively, we see that $A\ale\C_A^m(\bar H)$.\qed

\satz\label{fitting} Let $G$ be a type-definable group in a simple theory. Then
the Fitting subgroup $F(G)$ is nilpotent.\esatz
\bew As $F(G)$ is soluble by Lemma \ref{FGsoluble}, by Lemma \ref{soluble} there
is a chain
$$\{1\}=S_0<S_1<\cdots<S_d=S$$
of relatively definable normal subgroups of $G$ such that $S$ contains $F(G)$
and all quotients
$S_{i+1}/S_i$ for $i<d$ are abelian. Since $F(S_i)=F(G)\cap S_i$, we may assume
by induction on $d$ that $F(G)'\le F(S_{d-1})$ is nilpotent.
By Lemma \ref{nilpotent} and Remark \ref{nilpotentnormal} there is a relatively
definable normal nilpotent group $N$ containing $F(G)'$, and a relatively
definable series
$$\{1\}=N_0<N_1<\cdots<N_k=N$$
of normal subgroups of $G$ with $[N,N_{i+1}]\le N_i$ for all $i<k$.

Fix $i>0$. Any $g\in F(G)$ is contained in a normal nilpotent subgroup
$H_g$. Since $N_iH_g$ is again nilpotent, there is $m_g<\omega$
such that
$$(g-1)^{m_g}N_i\le N_{i-1}.$$
By Fact \ref{dcc} there is a finite tuple
$\bar g\in F(G)$ such that for any $g\in F(G)$ the index $|C_{N_i}(\bar
g/N_{i-1}):C_{N_i}(\bar g,g/N_{i-1})|$ is finite. Furthermore, by Lemma \ref{3}
(applied to $G/N$ with its abelian subgroup $F(G)/N$ acting on
$N_i/N_{i-1}$ by conjugation) there is
$m_i<\omega$ and a relatively definable group $H_i\ge F(G)$ such
that $N_i\ale\C_G^{m_i}(H_i/N_{i-1})$. Then, the finite intersection
$\bigcap_iH_i$ is a relatively
definable supergroup of $F(G)$. By Facts \ref{dcc} and \ref{sch} there is a
relatively definable normal subgroup $H$ which is a finite extension of a
finite intersection of $G$-conjugates of $\bigcap_iH_i$. Thus $H\ge F(G)$,
and $H\ale H_i$ for all $i$.
By Lemma \ref{soluble} applied to the abelian normal subgroup $F(G)/N$ of
$G/N$, we may restrict $H$ and assume that there are relatively definable
normal subgroups
$$N\le Z\le A\le H$$
of $G$ with $Z/N$ and $H/A$ finite and $A/Z$ abelian. Then
$$N_i\ale \C_G^{m_i}(H_i/N_{i-1})\le \C_G^{m_i}(H_i\cap
H/N_{i-1})=\C_G^{m_i}(H/N_{i-1}),$$
and inductively
$$\begin{aligned}N_i&\ale\C_G^{m_i}(H/N_{i-1}) \\ &\le
\C_G^{m_i}(H/\C_G^{m_{i-1}}(H/N_{i-2}))=\C_G^{m_i+m_{i-1}}(H/N_{i-2})\\
&\vdots\\
&\le\C_G^{m_i+m_{i-1}+\cdots+m_1}(H).\end{aligned}$$
Thus $N\ale\C_N^m(H)$ for $m=m_1+m_2+\cdots+m_k$. Since $Z\ale N$ and
$N\le A\le H$ we obtain $Z\ale\C_A^m(A)$, whence $A=\Z^{m+1}(A)$. So $A$ is
nilpotent-by-finite by Lemma \ref{nilpotent}, as is $F(G)$, since $F(G)\ale A$.
If $F$ is a normal nilpotent subgroup of finite index in $F(G)$ and $K$ a normal
nilpotent group
containing representatives for $F(G)/F$, then $FK$ is a nilpotent, as is $F(G)$.
\qed

\section{Hyperdefinable groups}
In the previous section we have systematically used the fact that
type-definability is preserved under quotients whenever we divide out by a
relatively definable subgroup. However, type-definability is not
preserved when quotienting by a type-definable subgroup, and in
fact such quotients (and even the slightly more general ones
defined below) arise naturally from model-theoretic considerations
in simplicity theory. We are thus led to the following definition.
\defn A {\em hyperdefinable} group is a group whose domain is given by a
partial type $\pi$ modulo a type-definable equivalence relation $E$, and whose
group law is induced by an $E$-invariant type-definable relation on
$\pi^3$.\edefn
Note that a quotient of a hyperdefinable group by a hyperdefinable group is
again hyperdefinable. In this context, we have
to replace {\em finite index} by {\em bounded index}, i.e.\ the index remains
bounded even in a very saturated elementary extension. With this replacement in
the definitions of Section \ref{sec1}, almost containment is transitive,
commensurability is an equivalence relation, and we still have
good definability properties in a simple theory.

\tats[{\cite[Proposition 4.4.10 and Corollary 4.5.16]{wa00}}] Let $K$ and $H$ be
hyperdefinable subgroups of a hyperdefinable group $G$ in a simple theory, with
$H\ale\N_G(K)$.
Then:\begin{enumerate}
\item $\N_G(K)$ is hyperdefinable, and a conjugacy-connected component $\tilde
K$ exists.
\item $\C_G^n(H/K)$ is hyperdefinable for all $n<\omega$.
\end{enumerate}\etats
Moreover, the proof of Proposition \ref{sym} remains valid in the
hyperdefinable context.

In contrast to the type-definable case, simplicity does not necessarly yield
a finite chain condition on centralizers (even though there is an ordinal
$\alpha$ such that any descending chain of hyperdefinable subgroups having
unbounded index in its predecessor stabilizes, up to bounded index, after
$\alpha$ many steps). In order to adapt the arguments from the
previous section we shall make a stronger assumption, supersimplicity. More
precisely, we shall assume the following consequence of supersimplicity:
There is no infinite descending chain of hyperdefinable subgroups, each of
unbounded index in its predecessor. In particular, we obtain a minimal condition
on centralizers, up to bounded index.

As a consequence, all proofs of the previous section adapt to this wider context
and therefore we obtain the same result, up to bounded index. Note that Remark
\ref{nilpotentnormal} need no longer hold, as a system of representatives for a
subgroup of bounded index can now be infinite.

Alternatively, we offer a distinct proof of virtual nilpotency of the
Fitting subgroup of a hyperdefinable group of ordinal $\SU$-rank in a simple
theory,
which in addition provides a bound on the nilpotency class. For the rest of the
section, the ambient theory will be simple. We first recall some
facts starting with the Lascar inequalities for $\SU$-rank.

\begin{fact}[{\cite[Theorem 5.1.6 (1)]{wa00}}]
If $H$ and $K$ are hyperdefinable subgroups of a common hyperdefinable group,
then
$$\SU(K)+\SU(HK/K)\le \SU(HK)\le \SU(K)\oplus \SU(HK/K),$$
where $\oplus$ is the least symmetric increasing function on ordinals satisfying
$f(\alpha,\beta+1)=f(\alpha,\beta)+1$.
\end{fact}

\begin{fact}[{\cite[Proposition 5.4.3]{wa00}}]\label{FactNormal1}
Let $G$ be an $\emptyset$-hyperdefinable group of rank
$\SU(G)=\omega^\alpha\cdot n+\gamma$ with $\gamma<\omega^\alpha$. Then $G$
has
an $\emptyset$-hyperdefinable normal subgroup $H$ of $\SU$-rank
$\omega^\alpha\cdot n$.
\end{fact}

\kor\label{Lemma:series}
Let $G$ be an $\emptyset$-hyperdefinable group of rank
$\SU(G)=\omega^{\alpha_1}\cdot n_1+\ldots+\omega^{\alpha_k}\cdot n_k$ with
$\alpha_i>\alpha_{i+1}$ for $i<k$ and $n_i>0$ for $i\le k$. Then there exists a
series of $\emptyset$-hyperdefinable $G$-invariant subgroups $$\{1\}=G_0\unlhd
G_1 \unlhd \cdots\unlhd G_\ell=G$$ with $\ell\le n_1+\ldots + n_k$ such that
each quotient
$G_{i+1}/G_i$ is unbounded of monomial $SU$-rank $\omega^{\beta_i}\cdot m_i$
and its $\emptyset$-hyperdefinable $G$-invariant subgroups of unbounded index
have $SU$-rank strictly smaller than $\omega^{\beta_i}$.
\ekor
\bew By Fact \ref{FactNormal1} there is an $\emptyset$-hyperdefinable normal
subgroup $G_1$ of $G$ of minimal monomial Lascar rank of the form
$\SU(G_1)=\omega^{\alpha_1}\cdot m$ with positive $m\le n_1$. By
minimality, $G_1$ is as required. If $\SU(G_1)=\SU(G)$ we are done.
Otherwise, $\SU(G/G_1)<\SU(G)$ by the Lascar inequalities, so we finish by
induction on $\SU(G)$.\qed

Next, we recall the supersimple version of Zilber's Indecomposability Theorem.

\begin{fact}[{\cite[Theorem 5.4.5 and Remark 5.4.7]{wa00}}]\label{FactZilber}
Let $G$ be an $\emptyset$-hyper\-definable group of rank
$SU(G)<\omega^{\alpha+1}$. If $\mathfrak X$ is a family of
hyperdefinable subsets of $G$, then there exists a hyperdefinable subgroup
$K\le X_1^{\pm 1}\cdots X_m^{\pm 1}$ for some $X_1,\ldots,X_m\in \mathfrak X$,
such that $\SU(XK)<\SU(K)+\omega^\alpha$ for all $X\in\mathfrak X$. Moreover,
$\SU(K)=\omega^\alpha\cdot n$, and $K$ is unique up to commensurability. In
particular, if $\mathfrak X$ is invariant under all automorphisms we can choose
$K$ hyperdefinable over $\emptyset$, and if $\mathfrak
X$ is $G$-invariant, we can take $K$ to be normal in $G$.\end{fact}

Finally, we state the hyperdefinable version of our main result in the
supersimple case. Recall that the $\emptyset$-connected component
$N_\emptyset^0$ of a hyperdefinable group $N$ is the intersection of all
$\emptyset$-hyperdefinable subgroups of $N$ of bounded index in~$N$.

\satz\label{hyperfitting} Let $G$ be an $\emptyset$-hyperdefinable group of rank
$\SU(G)=\omega^{\alpha_1}\cdot n_1+\ldots+\omega^{\alpha_k}\cdot n_k$. Then
$F(G)$ has bounded index in an $\emptyset$-hyperdefinable FC-nilpotent normal
subgroup $N$ of
class $\ell\le n_1+\ldots+n_k$. In particular, the
hyperdefinable normal group $N_\emptyset^0$ is nilpotent of class $2\ell$ and
has bounded index in $F(G)$.
\esatz
\bew By Lemma \ref{Lemma:series} there is a finite series of
$\emptyset$-hyperdefinable $G$-invariant subgroups
$$\{1\}=G_0 \unlhd G_1 \unlhd \cdots \unlhd G_\ell=G$$ with $\ell\le
n_1+\cdots+n_k$, such that each quotient $G_{i+1}/G_i$ is unbounded of monomial
Lascar
rank $\omega^{\beta_i}\cdot m_i$, and its $\emptyset$-hyperdefinable
$G$-invariant subgroups of unbounded index have $\SU$-rank strictly smaller than
$\omega^{\beta_i}$. Clearly, we may assume that all $G_i$ are
$\emptyset$-connected, i.e.\ have no $\emptyset$-hyperdefinable subgroup of
bounded index.

Let $N$ be the intersection $\bigcap_{i<\ell}\C_G(G_{i+1}/G_i)$, an
$\emptyset$-hyperdefinable normal subgroup of $G$.
Note that $N\le\C_G((G_{i+1}\cap N)/(G_i\cap N))$. Hence by symmetry we get
$$G_{i+1}\cap N\ale\C_G(N/(G_i\cap N))$$
for all $i<\ell$. Inductively,
$$\begin{aligned}N&=G_\ell\cap N\ale\C_G(N/(G_{\ell-1}\cap N))\\
&\le\C_G(N/(\C_G(N/(G_{\ell-2}\cap N)))=\C_G^2(N/(G_{\ell-2}\cap N))\\
&\le\C_G^2(N/\C_G(N/(G_{\ell-3}\cap N)))=\C_G^3(N/(G_{\ell-3}\cap N))\\
&\le\cdots\le\C_G^\ell(N/(G_0\cap N))=\C_G^\ell(N).\end{aligned}$$
It follows that $N=\Z_\ell(N)$ is FC-nilpotent of class at most $\ell$. Then
$N_\emptyset^0$ is nilpotent of class $2\ell$ by \cite[Proposition 4.4.10
(3)]{wa00}.

In order to finish, we shall show $F(G)\le N$. Fix $i\le\ell$,
and consider the $\emptyset$-invariant family
$\mathfrak X$ formed by the hyperdefinable sets $X_a=[a,G_{i+1}]/G_i$ for $a\in
F(G)$. Suppose, towards a contradiction, that the $\SU$-rank of some of these
sets is greater than
$\omega^{\beta_i}$. By Fact \ref{FactZilber} applied to $G_{i+1}/G_i$
we obtain an $\emptyset$-hyperdefinable $G$-invariant subgroup $H\le
G_{i+1}/G_i$ of monomial $\SU$-rank which is contained in a finite product of
sets $X_{a_0}^{\pm1},\ldots, X_{a_m}^{\pm1}$
from $\mathfrak X$; moreover, $\SU(H)\ge\omega^{\beta_i}$. Thus $H$ has bounded
index in
$G_{i+1}/G_i$ and must be equal by $\emptyset$-connectivity. But every
$a_jG_i$ is contained in a normal nilpotent subgroup of $G/G_i$
which must also contain
$X_{a_j}$, so $K=\langle a_jG_i,X_{a_j}:j\le m\rangle$ is a nilpotent subgroup
of $G/G_i$. However,
\begin{equation*}\tag{\dag}H\le X_{a_0}^{\pm 1}\cdots X_{a_m}^{\pm
1}=[a_0,H]^{\pm1}\cdots[a_m,H]^{\pm1}\subseteq H\end{equation*}
and we must have equality, contradicting nilpotency of $K$: If $H$ is in the
$k$-th element $\gamma_k(K)$ of the lower central series of $K$, then equation
$(\dag)$ implies that $H\le\gamma_{k+1}(K)$. Thus
$H\le\bigcap_{k<\omega}\gamma_k(K)=\{1\}$, a contradiction.

It follows that $\SU(X_a)<\omega^{\beta_i}$ for all $X_a\in\mathfrak X$. As
$X_a$ is in bijection with $G_{i+1}/C_{G_{i+1}}(a/G_i)$, the Lascar inequalities
imply that $X_a$ is bounded and so $a\in\C_G(G_{i+1}/G_i)$. Thus
$F(G)\le N$, as required.\qed

\end{document}